\documentclass[letterpaper, 12pt]{article}

\usepackage[margin=1in]{geometry}
\usepackage{setspace}\onehalfspace
\usepackage{microtype}

\usepackage{amsmath} 
\allowdisplaybreaks
\usepackage{amssymb}
\usepackage{amsthm}
\usepackage{algpseudocode}
\usepackage{algorithm}
\usepackage{cite}
\usepackage{multirow}

\theoremstyle{plain}
\newtheorem{thm}{Theorem}[section]

\numberwithin{equation}{section} \theoremstyle{definition}

\newtheorem{alg}[thm]{Algorithm}

\newtheorem{exmp}[thm]{Example}

\numberwithin{equation}{section}

\usepackage{booktabs} 
\usepackage{graphicx}
\usepackage[counterclockwise]{rotating} 
\usepackage{subcaption} 
\usepackage[normalem]{ulem} 
\usepackage{showexpl}
\usepackage{threeparttable} 
\usepackage{lipsum}
\usepackage[hyphens]{url}
\usepackage[colorlinks,urlcolor=blue]{hyperref}

\usepackage{pgf}
\pgfdeclarelayer{background}
\pgfdeclarelayer{foreground}
\pgfsetlayers{background,main,foreground}

\usepackage{authblk}

\newcommand{\email}[1]{{\href{mailto:#1}{\nolinkurl{#1}}}}

\usepackage{bm}

\usepackage{xcolor}

\title{A Practical Approach to Quasi-convex Optimization}
\author[2]{Sompong Dhompongsa}
\author[1,2]{Poom Kumam\footnote{Corresponding Author: \email{poom.kum@kmutt.ac.th }}}
\affil[1]{Fixed Point Research Laboratory,  Fixed Point Theory and Applications Research Group, Center of Excellence in Theoretical and Computational Science (TaCS-CoE), Faculty of Science, King Mongkut’s University of Technology Thonburi (KMUTT), 126 Pracha Uthit Rd., Bang Mod, Thung Khru, Bangkok 10140, Thailand}
\affil[2]{Center of Excellence in Theoretical and Computational Science (TaCS-CoE), Faculty of Science, King Mongkut’s University of Technology Thonburi (KMUTT), 126 Pracha Uthit Rd., Bang Mod, Thung Khru, Bangkok 10140, Thailand}

\date{}

\begin{document}
\maketitle
\vspace{-1cm}
\begin{abstract}
A new and simple method for quasi-convex optimization is introduced from which its   various applications can be derived. Especially, a global optimum under constrains can be approximated for all continuous functions.
\end{abstract}

\section{Problems and Solutions}
A classical form of optimization formulated as: To find
\begin{equation} \label{eq main}
\begin{aligned}
& \underset{x}{\arg \min}f(x) \\
\text{subject to } & g_i(x) = 0, \\
& h_j(x) \leq 0, \quad (i=1,\dots,m;  j=1,\dots,n).\\
\end{aligned}
\end{equation}
Following \cite{dhompongsa2020practical}, a solution of \eqref{eq main} can be obtained by minimizing a deformed function
\begin{equation}\label{eq 2}
	f_t (\cdot | K,M) := (1-t)(f(\cdot)-K)+tM \left [  \sum_{i=1}^{m}  |g(\cdot)| + \sum_{j=1}^{n}(h (\cdot) + |h (\cdot)|) \right], \quad t \in (0,1)
\end{equation}
$t$ is very close to $1$ and $K,M$ are large. A problem, however, arrives when the function
$$F(\cdot) : = \left [  \sum_{i=1}^{m}  |g(\cdot)| + \sum_{j=1}^{n}(h (\cdot) - |h (\cdot)|) \right],$$
is not convex. A minimizer may be too far away from the feasible set. To fix the problem, we replace $f_t$ by $f_t-|f_t |.$ The rest of the paper will now pay attention to look for a method for optimization of quasi-convex like functions.

\begin{figure}[h!]
	\begin{subfigure}{\textwidth}
		\centering
		\includegraphics[width=1\linewidth]{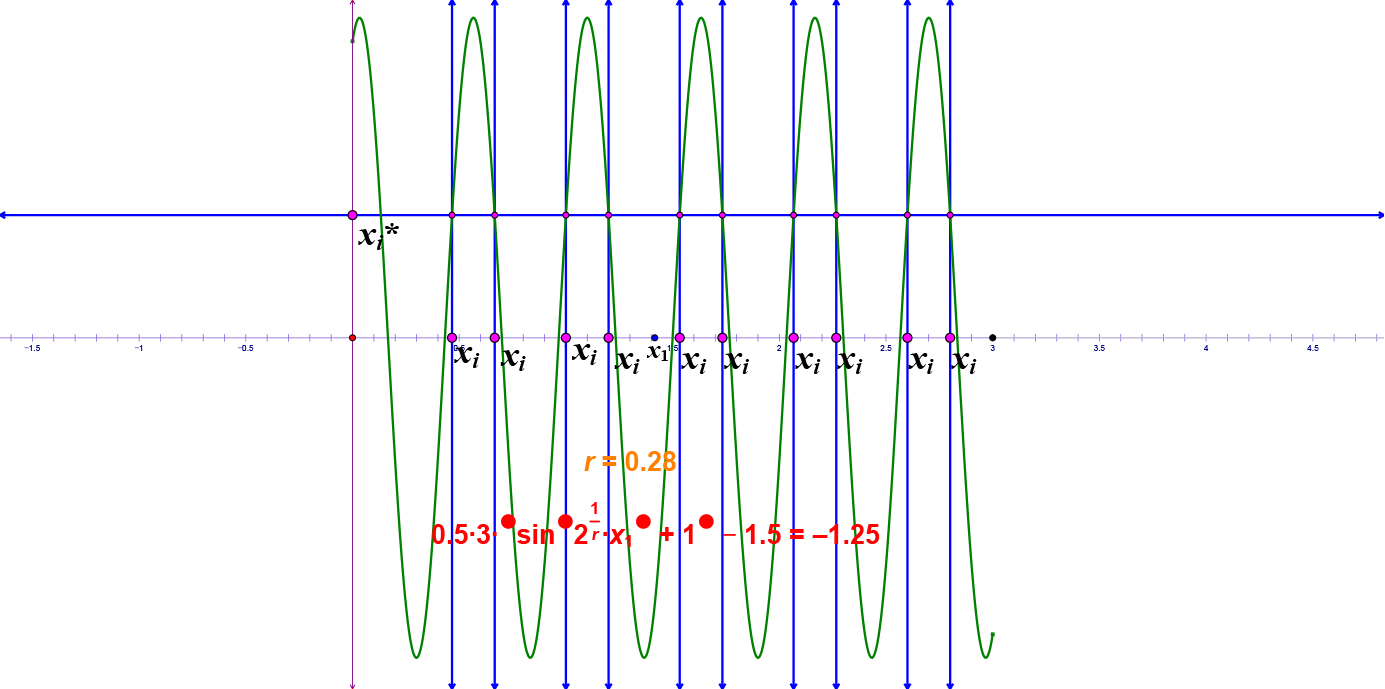}
			\caption{}
				\label{Fig:1(a)}
	\end{subfigure}\\
	\begin{subfigure}{\textwidth}
		\centering
		\includegraphics[width=1\linewidth]{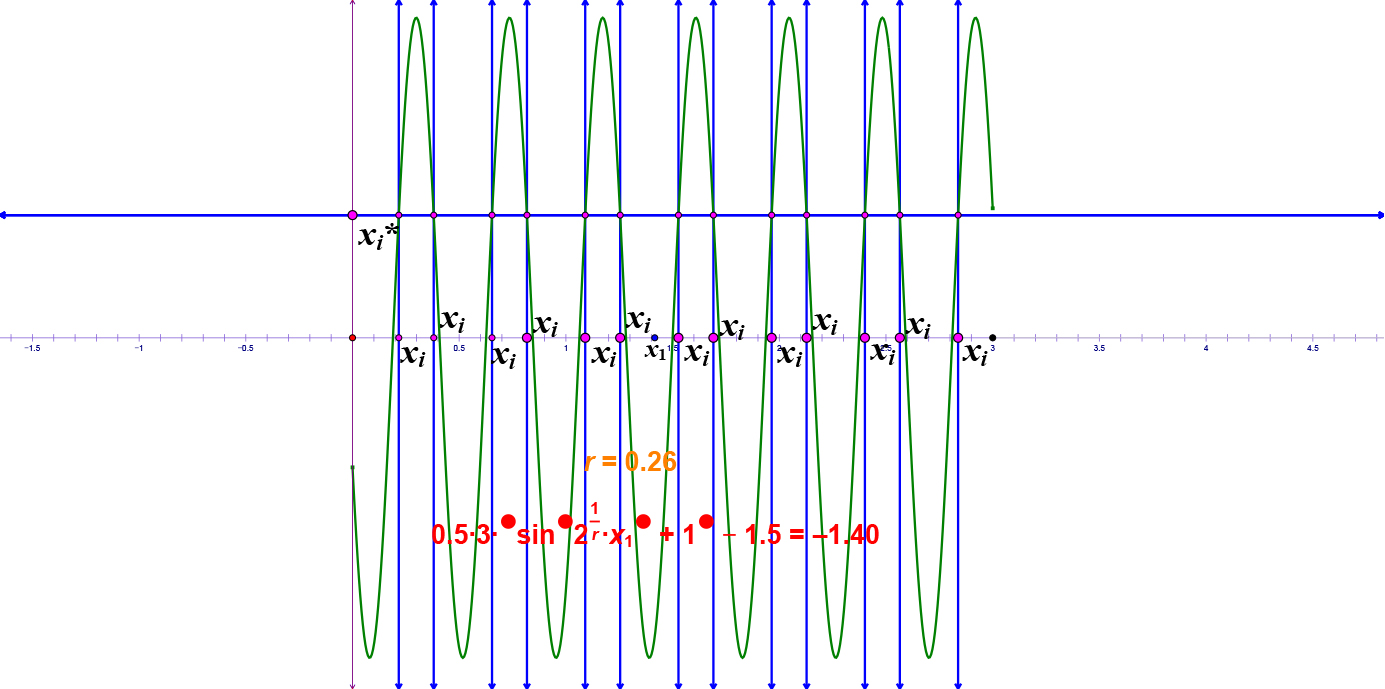}
			\caption{}
				\label{Fig:1b}
	\end{subfigure}
	\caption{}
		\label{Fig1}
\end{figure}

Consider a function $g: C = \prod_{i=1}^{p}[a_i,b_i] \to (-\infty,0]$.  Suppose $g<0$ on a small neighborhood containing $x^*=(x_1^*,\ldots ,x_p^*).$ To find a point in this neighborhood we introduce a point $x^{'}=(x_1^{'}, \ldots ,x_p^{'})$ which transforms a point $x=(x_1,\ldots,x_p) \in C$ under the rule:      
\begin{equation}\label{eq3}
	x_i^{'} = .5(\sin(2^{\frac{1}{r_i}}x_i)+1)(b_i - a_i)+a_i, \quad    i=1, \ldots, p,
\end{equation}                                
where $r=(r_1,\ldots,r_p) \in (0,\infty)^p.$ \\
It is seen that $x_i^{'}= \gamma(\beta(x_i))$ where $\gamma(\beta)$ is the composition of $\beta : [a_i, b_i] \to [-1,1], b  \mapsto \sin (2^{\frac{1}{r_i}})b)$  followed by $ \gamma : [-1,1] \to [a_i,b_i ], c \mapsto .5(b_i-a_i )(c+1)+a_i. $ Observe that $\gamma(\beta)$ sends $[a_i,b_i ]$ onto itself many to one.   \\
Put $u(x,r)=g(x^{'} ).$ Clearly, if  $x_i^{'}= x_i^{*}$, then $u(x,r)<0$ for all $x=(x_1,\ldots,x_p)$ for which the
 relation $x_i^{'}= x_i^{*}$ holds for each $i$ 
(see Figure \ref{Fig1}). Figure  \ref{Fig:1(a)} and \ref{Fig:1b} also reveal the fact that the number of points satisfying the relation $x_i^{'}= x_i^{*}$, increases as $r_i$ gets smaller.

\newpage
\noindent This shows that, instead of working on $f$, we work on $u(x,r)$ to get more chances in finding a solution, i.e., a point in the neighborhood of $x^{*}$. Still, the problem remains, since we have no clue to move from any initial point $(x_0,r_0)$= $(x_{01},\ldots,x_{0p},r_{01},\ldots,r_{0p})$ where $u(x_0,r_0)=0$	 to a new point $(x_1,r_1)$ with $u(x_1,r_1)<0$ (see Figure \ref{Fig:4}). Fortunately, we only add to $u(x,r)$  with a convex function
$$v(x,r)=(c_1-x_1)^2 + \cdots + (c_p-x_p )^2+(2-r_1 )^2+ \cdots+(2-r_p)^2$$
and call $w=u+v$. Here $c_i$ is the mid-point of $[a_i,b_i]$ for all $i$ and $2$ can be replaced by any positive number.  Finally, we can apply any convex optimization algorithm to get any (local) minimizer of $w$ over $(x_1,\ldots,x_p,r_1,\ldots,r_p).$ Note that any minimizer $x=(x_1,\ldots,x_p)$ solves the equation $x_i^{'}= x_i^{*}$ at each $i$.  \\
Now instead of finding $\underset{x}{\arg \min}f(x)$, we set to find $x^{*}$ for which $f(x^{*} )<k$ for any given $k.$ Then put $g=f-k-|f-k|.$  So we shall find $x^{*}$ satisfying  $g(x^{*} )<0.$ This function $g$ includes function $f_t$ mentioned above.  \\
In the following computation, we set $r_1= \cdots =r_p.$  Moreover, since $(c_1-x_1)^2+\cdots+(c_p-x_p )^2$  is convex as a function of $x$ with its minimum value occurs at the point $c=(c_1,\ldots,c_p),$ we may choose to concentrate on a piece of the graph of w around $(c_1,\ldots,c_p).$  
\begin{alg} (To find $x^{*}$for which $f(x^{*} )<k$) \\
	Given $f : C = \prod_{i=1}^p [a_i,b_i] \to \mathbb{R},  k \in \mathbb{R}$. Set $g=f-k-|f-k|.$  \\
	Plug in $x{'}$ by \eqref{eq3} to get $g(x^{'})$. Recall that we fix all $r_1=\ldots=r_p$ and set  $v(x)= (c_1-x_1)^2+\cdots+(c_p-x_p )^2+(2-r_1)^2,$ then compute  $w(x,r_1 ):= g(x^{'} )+v(x,r_1 ).
$ \\
	Minimize $w$ over $(x,r_1 )$ setting $(a_1,\ldots,a_p,r_1)$ as our initial point where $r_1$ is close to $0$ and get any  minimizer $(x^{*},r_1^{*}).$ We then get $g(x^{*'}) )<0,$ i.e., $f(x^{*'}) <k.$
\end{alg}
\noindent Obviously, a global minimum of $f$ can be obtained by decreasing $k$ as much as we can. \\
Another way to minimize w is as follow: We are looking for cross sections of w at the point $(c,r_1^{*} )$ where we find first a minimizer $r_1^{*}$ of the function $r_1 \mapsto  w(c,r_1 ).$ This means that the cross section along $r_1$ where the variables $x=(x_1,\ldots,x_p)$ is fixed at $c=(c_1,\ldots, c_p)$ has a local minimum at $r_1^{*}.$ The point  $(c,r_1^* )$ will automatically a minimizer of $w$. Actually, we can set any point to start with instead of starting at $c$.

\begin{exmp}
	Note, in Figure \ref{Fig:2}, \ref{Fig:3}, \ref{Fig:4}, \ref{Fig:5a}, \ref{Fig:5b}, that $(x,y,z,u)$  stands for $(x_1,\ldots x_4 )$ and set $r_1=r_2=r_3=r_4=r$ (don’t be mixed up this $r$ and the vector $r=(r_1,r_2,r_3,r_4)$ in the above context). \\
	Consider  the function \\
	$f(x_1,\ldots, x_4 ) = (x_1-2)^2+(x_2+2)^2+(x_3-2)^2+(x_4+2)^2-.55-|(x_1-2)^2+(x_2+2)^2+(x_3-2)^2+(x_4+2)^2-.55|,  ((x_1,\ldots, x_4 )) \in [-1.5,1.5]^4 $
	(see Figure \ref{Fig:2} for its graph). Find $((x_1,\ldots, x_4 ))  \in [-1.5,1.5]^4$ satisfying $f((x_1,\ldots, x_4 ))<0.$
		\begin{figure}[h!]
			\centering
			\includegraphics[width=1\linewidth]{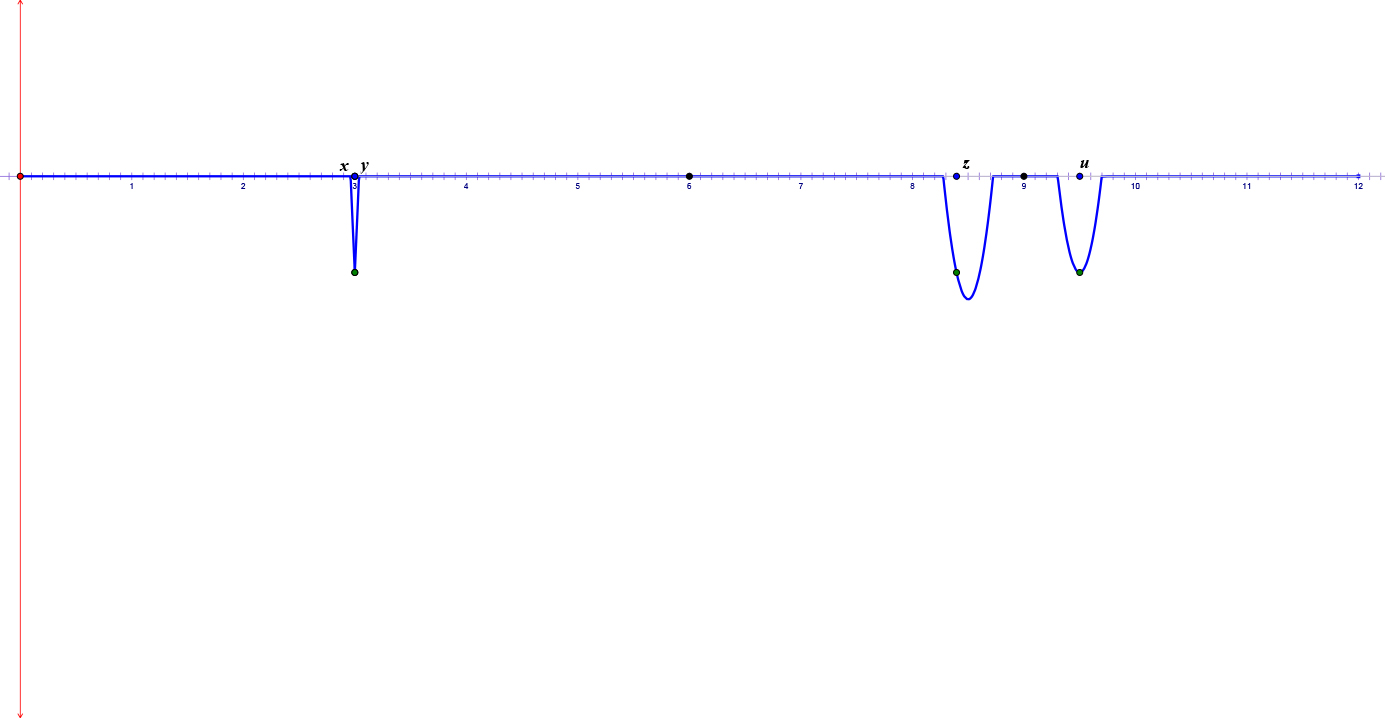}
			\caption{}
			\label{Fig:2}
	\end{figure}
	It is easy to guess that $(1.5,-1.5,1,-1)$ is one such an answer. 
	Most of the time, all cross sections at arbitrary $(x_1,\ldots,x_4 )$ of $f$ will be zero as in Figure \ref{Fig:3}.
			\begin{figure}[h!]
		\centering
		\includegraphics[width=1\linewidth]{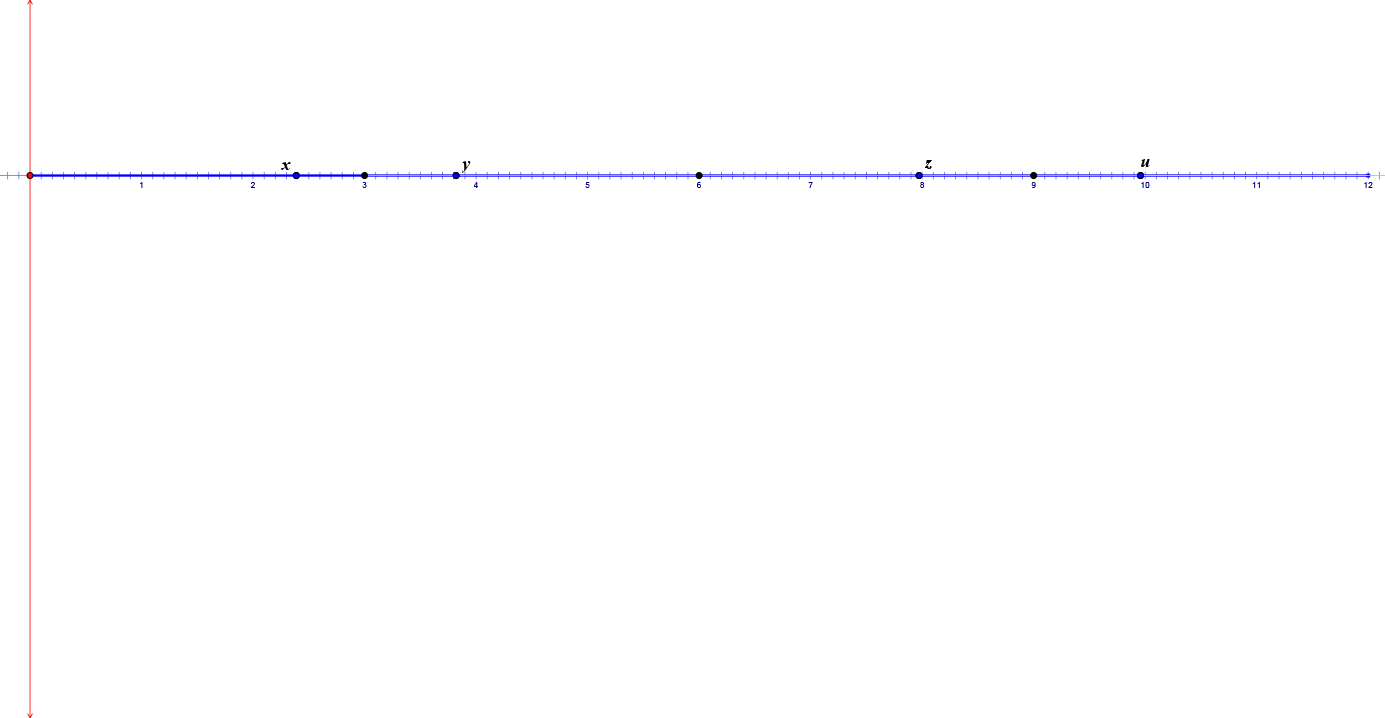}
		\caption{}
		\label{Fig:3}
	\end{figure}

	Plug in $(x_1^{'},\ldots,x_4^{'} )$ by \eqref{eq main} into f and get its graph as in Figure \ref{Fig:4}:
	\newpage
	\clearpage
				\begin{figure}[h!]
		\centering
		\includegraphics[width=1\linewidth]{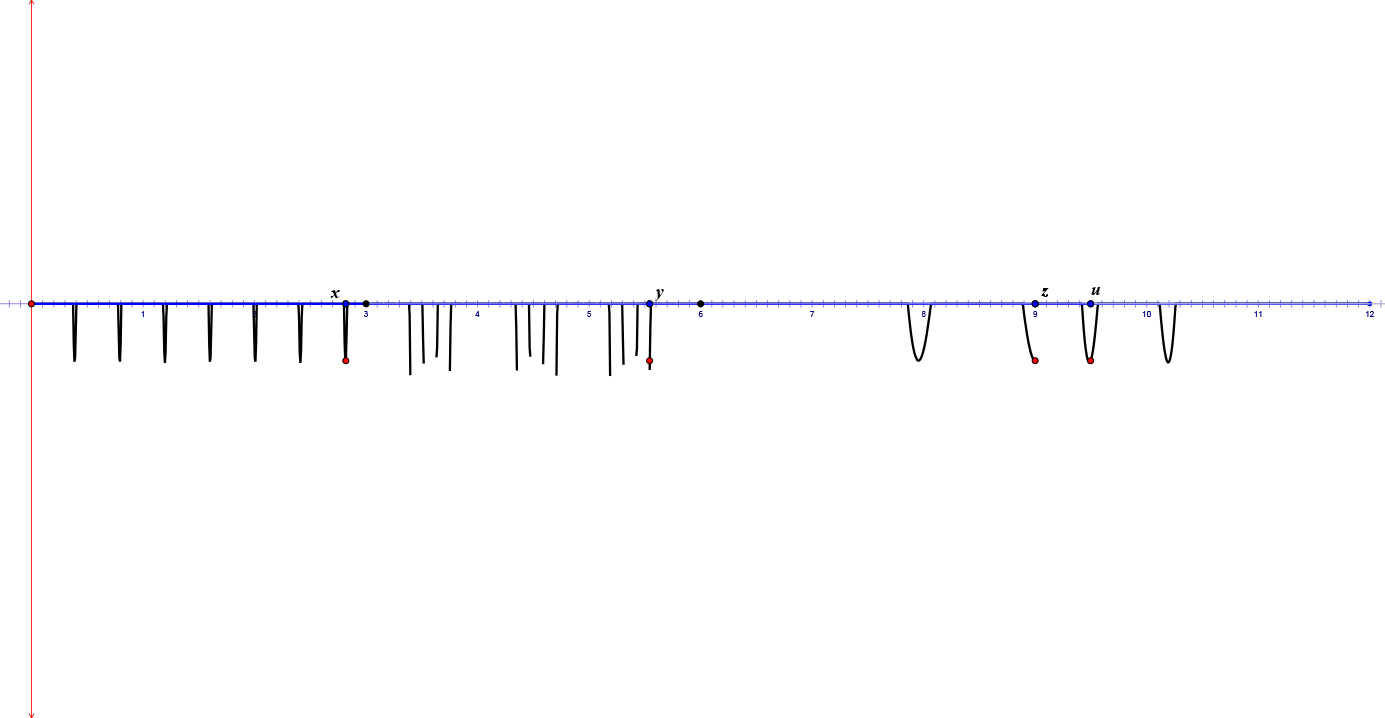}
		\caption{}
		\label{Fig:4}
	\end{figure}
\noindent Now consider the graphs of $w$ in Figure \ref{Fig:5a} and Figure \ref{Fig:5b} for different minimizers $r_1^{*}$. Each Figure
 shows a point in a neighborhood of $x^{*}$ where $f(x^{*} )<0$, namely, the point $(1.5,-1.48,1,-1)$ and the point $(1.6,-1.5,.89.-1)$ respectively.  
	\begin{figure}[h]
	\begin{subfigure}{\textwidth}
		\centering
		\includegraphics[width=1\linewidth]{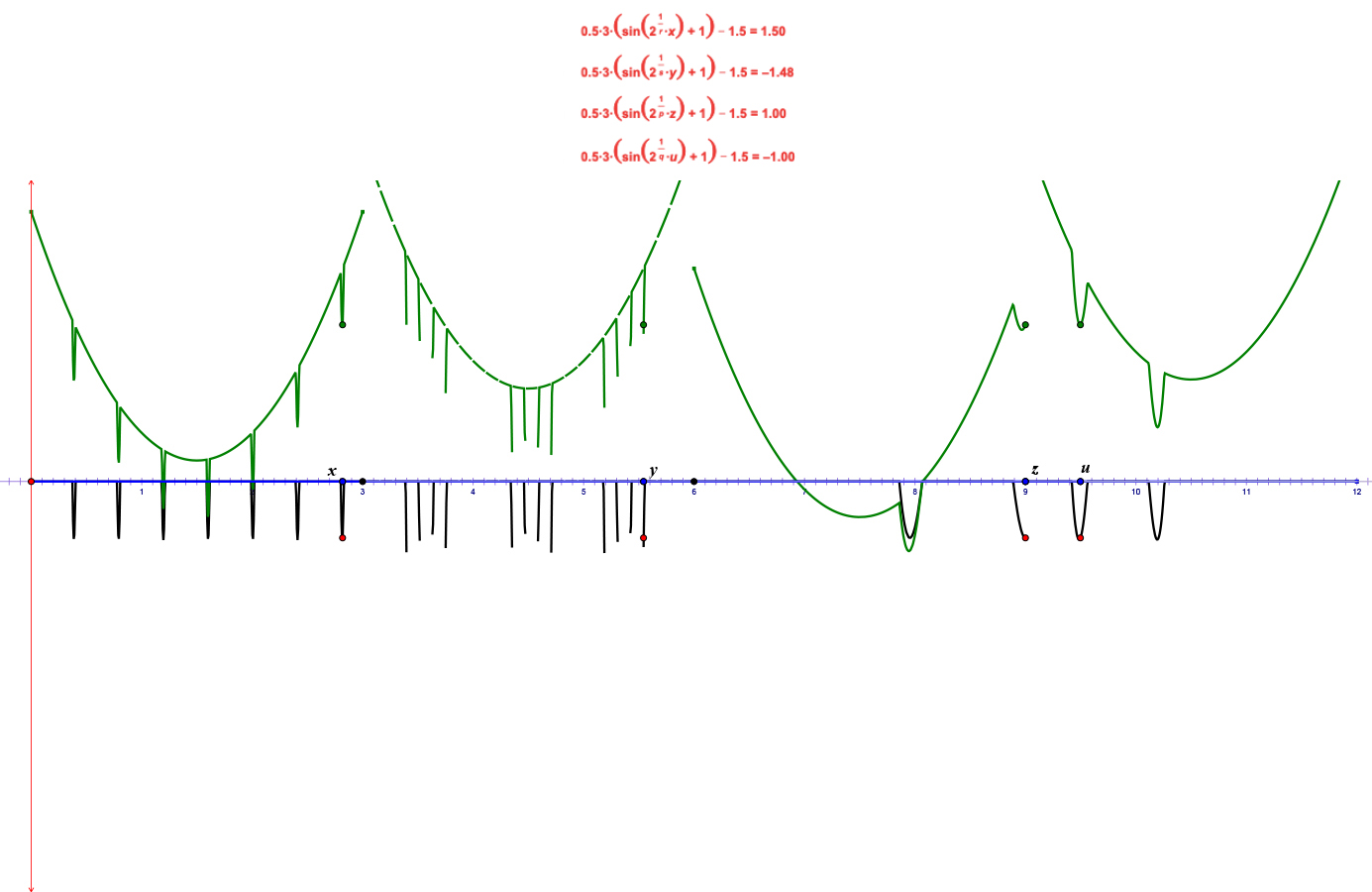}
		\caption{}
		\label{Fig:5a}
	\end{subfigure}\\
	\begin{subfigure}{\textwidth}
		\centering
		\includegraphics[width=1\linewidth]{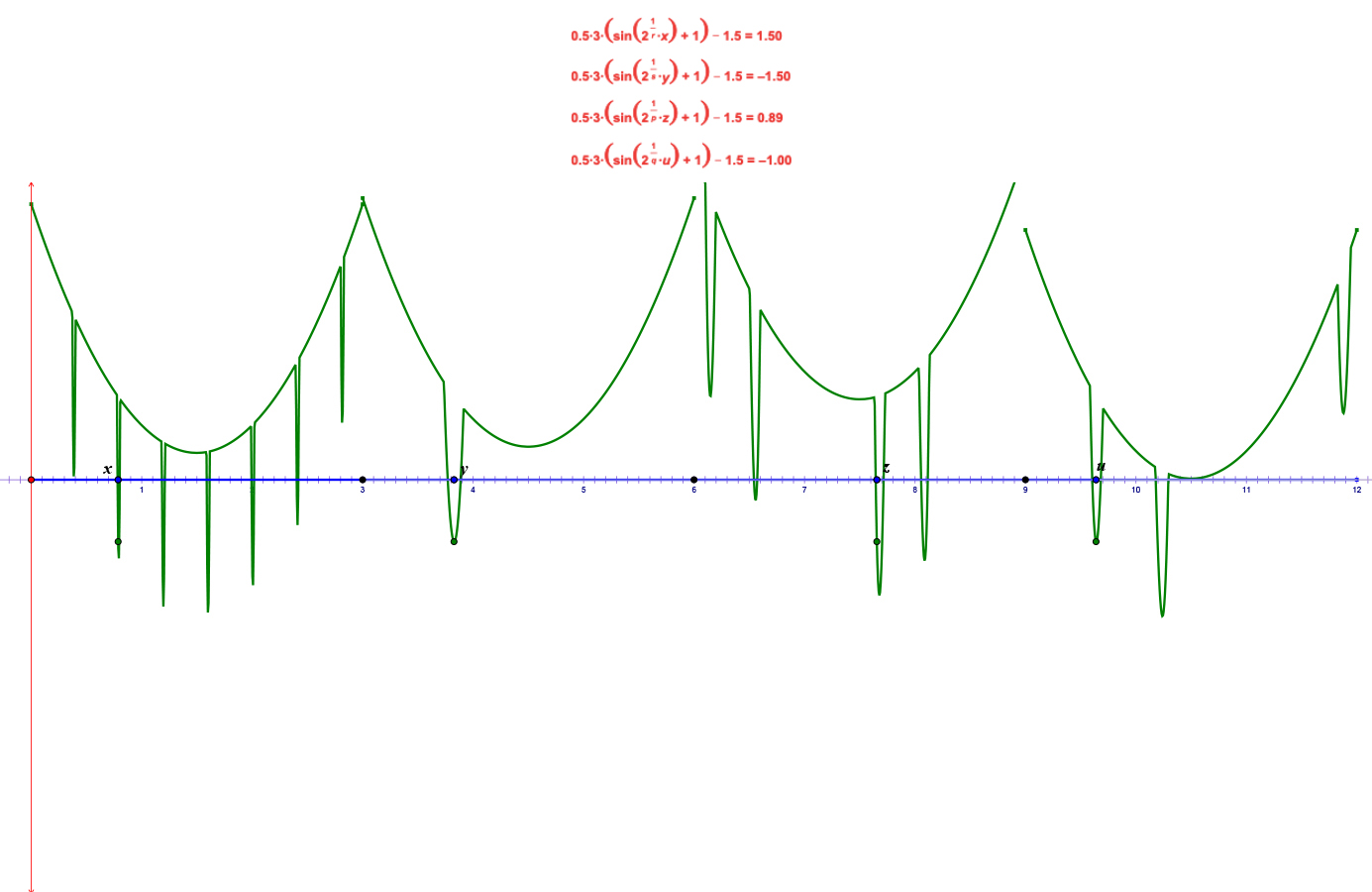}
		\caption{}
		\label{Fig:5b}
	\end{subfigure}
	\caption{}
	\label{Fig:5}
\end{figure}
\end{exmp}

\newpage
\clearpage
\section*{Acknowledgments}
 The authors acknowledge the financial support provided by the Center of Excellence in Theoretical and Computational Science (TaCS-CoE), KMUTT. 
Moreover, this research project is supported by Thailand Science Research and Innovation (TSRI) Basic Research Fund: Fiscal year 2021 under project number 64A306000005.
\bibliographystyle{vancouver}
\bibliography{ref}

\begin{thebibliography}{1}

\bibitem{dhompongsa2020practical}
Dhompongsa S, Jirakitpuwapat W, Khammahawong K, Kumam P.
\newblock A practical approach to optimization.
\newblock arXiv preprint arXiv:200910452. 2020;.

\end{thebibliography}
\end{document}